\documentclass[12pt,a4paper]{amsart}
\usepackage{amsmath,amssymb}
\usepackage[T1]{fontenc}
\usepackage[utf8]{inputenc}
\usepackage{lmodern}
\usepackage[numbers,sort&compress]{natbib}
\usepackage{hyperref}
\hypersetup{hypertexnames=false,hidelinks}
\newcommand{\T}{\mathbb{T}}
\newcommand{\N}{\mathbb{N}}
\newcommand{\ldimM}{\underline{\dim}_{M}}

\newtheorem{thm}{Theorem}[]
\newtheorem{cor}[thm]{Corollary}
\newtheorem{lem}[thm]{Lemma}
\newtheorem{prop}[thm]{Proposition}
\newtheorem{de}[thm]{Definition}
\theoremstyle{remark}
\newtheorem{rem}[thm]{Remark}
\date{}
\title[Rajchman measures as $H^1(\T^d)$-multipliers]{The $H^1(\T^d)$-multiplier norm of Rajchman measures, Zafran-type phenomena, and Walsh--Stieltjes series}
\author{Przemysław Ohrysko}
\thanks{Corresponding author: Przemysław Ohrysko.}
\address{University of Warsaw, Faculty of Mathematics, Informatics and Mechanics, Banacha 2, 02-097 Warsaw, Poland}
\email{p.ohrysko@gmail.com}
\author{Michał Wojciechowski}
\address{Institute of Mathematics, Polish Academy of Sciences, 00-956 Warsaw, Poland}
\email{miwoj.impan@gmail.com}
\author{Bartłomiej Zawalski}
\address{Department of Mathematics, Applied Mathematics, and Statistics, Case Western Reserve University, Cleveland, OH 44106, USA}
\email{bxz428@case.edu}
\begin{document}
\renewcommand{\shortauthors}{Ohrysko, Wojciechowski, and Zawalski}
\begin{abstract}
We study convolution by complex Rajchman measures as Fourier multipliers on the mean-zero part of the classical isotropic Hardy space $H^1(\T^d)$. Our main technical result gives a lower bound for the multiplier norm in terms of the amount of mass carried by sets of small lower Minkowski dimension. The estimate holds for every $d\ge1$; the argument is first established for real measures and then extended to complex measures by a finite phase decomposition. Since the testing functions have mean zero, a corresponding lower bound also holds on the full isotropic Hardy space. The estimate yields multidimensional Zafran-type criteria for non-natural spectra and, in particular, shows that positive Rajchman probability measures carried by sets of lower Minkowski dimension zero induce convolution operators with non-natural spectrum. In dimension one, for real measures, we also obtain the analogous estimate on the analytic space $H^1_0(\T)$ and frequency-localized failures of the natural-spectrum property. Finally, the binary form of the underlying combinatorial lemma gives Haar block estimates, which transfer to uncertainty principles and carrier-dimension results for Walsh--Stieltjes coefficients.
\end{abstract}

\maketitle

\section{Introduction}
The classical Wiener--Pitt phenomenon shows that the spectrum of a measure in the convolution algebra $M(\T)$ may be strictly larger than the closure of the range of its Fourier--Stieltjes transform \citep{wp}. Zafran developed a general spectral theory for multipliers whose transforms vanish at infinity \citep{Zafran} and, on the one-dimensional circle, proved that convolution by certain infinite Bernoulli convolutions satisfying suitable arithmetic conditions has non-natural spectrum on $H^1$ and on Lipschitz spaces \citep{Zaf2}. The purpose of the present paper is to replace these special one-dimensional constructions by a quantitative geometric mechanism valid on every torus $\T^d$: concentration of a Rajchman measure on sets of small lower Minkowski dimension forces a large convolution norm on $H^1(\T^d)$ and, through convolution powers, produces non-natural spectra. In particular, our zero-dimensional carrier criterion and the more general fractal-survival criterion below are genuinely multidimensional extensions of this circle of results.

Fix an integer $d\ge1$. Let $M(\T^d)$ denote the Banach algebra of all complex-valued finite Borel measures on the $d$-dimensional torus, equipped with the total variation norm and convolution. We write $M_0(\T^d)$ for its closed ideal of Rajchman measures, that is, measures $\mu$ for which $\widehat\mu(k)\to0$ as $|k|\to\infty$ in $\mathbb Z^d$.

Throughout, Haar measure on $\T^d=\mathbb R^d/\mathbb Z^d$ is normalized to have total mass one, and we use the Fourier convention
\begin{equation*}
 \widehat\mu(k)=\int_{\T^d}e^{-2\pi i k\cdot x}\,d\mu(x),
 \qquad k\in\mathbb Z^d.
\end{equation*}

Throughout the paper we use the lower Minkowski dimension. If $A\subset\T^d$ and $N(A,\delta)$ denotes the least number of cubes of side $\delta$ needed to cover $A$, we put
\begin{equation*}
 \ldimM A=\liminf_{\delta\downarrow0}
 \frac{\log N(A,\delta)}{\log(1/\delta)}.
\end{equation*}
Cubes may equivalently be replaced by balls. We shall also use that taking the closure of a set does not change its lower Minkowski dimension; a proof is recorded in Lemma \ref{lem:closure-minkowski}. For $0<\beta<d$ and $\mu\in M(\T^d)$ define
\begin{equation}\label{eq:cbeta-d}
 c^{(d)}_\beta(\mu)
 =\sup\bigl\{|\mu|(A):A\subset\T^d\text{ is Borel and }\ldimM A<\beta\bigr\}.
\end{equation}
Thus $c^{(d)}_\beta(\mu)$ measures the largest part of the mass of $\mu$ which can be detected on sets of lower Minkowski dimension below $\beta$. The lower, rather than upper, Minkowski dimension is the relevant quantity because the proof uses efficient dyadic coverings along infinitely many scales. When $d=1$ we omit the superscript and write $c_\beta(\mu)$.

Let $H^1_{\rm at,0}(\T^d)$ be the codimension-one, mean-zero part of the classical isotropic complex atomic Hardy space on $\T^d$. The real-variable Hardy-space framework underlying this formulation goes back to Fefferman and Stein \citep{FeffermanStein}; for the atomic description used here, see also \citep{Stein,SleddStegenga}. Thus it is the completion of finite sums $f=\sum_j\lambda_j a_j$, where the $a_j$ are complex atoms supported on cubes $Q_j\subset\T^d$ and satisfy
\begin{equation*}
 \int_{Q_j} a_j=0,
 \qquad \|a_j\|_\infty\le |Q_j|^{-1},
\end{equation*}
with norm given by the infimum of $\sum_j|\lambda_j|$ over all such decompositions. Its real subspace is the usual real-variable atomic Hardy space. All testing functions constructed below are real and are uniformly bounded multiples of atoms. Because they have mean zero, every lower estimate proved on $H^1_{\rm at,0}(\T^d)$ also gives the corresponding lower estimate on the full isotropic Hardy space, up to the harmless equivalence of the standard norms. We use the mean-zero space in the spectral statements in order to remove the zero Fourier mode. The main result is the following.

\begin{thm}[Main theorem]\label{dwa}\label{thm:multidim-main-intro}
Let $d\ge1$ and $0<\beta<d$. There is a constant $C(d,\beta)>0$ such that, for every complex measure $\mu\in M(\T^d)$,
\begin{equation}\label{eq:main-multidim}
 \|T_\mu:H^1_{\rm at,0}(\T^d)\longrightarrow L^1(\T^d)\|
 \ge C(d,\beta)c^{(d)}_\beta(\mu).
\end{equation}
Consequently, with another constant depending only on $d$ and $\beta$,
\begin{equation*}
 \|T_\mu:H^1_{\rm at,0}(\T^d)\longrightarrow H^1_{\rm at,0}(\T^d)\|
 \ge C(d,\beta)c^{(d)}_\beta(\mu).
\end{equation*}
Moreover, the lower bound is realized by uniformly bounded multiples of atoms at infinitely many dyadic scales.
\end{thm}

For real measures the proof is genuinely $d$-dimensional. A dyadic cube has $2^d$ children, and a non-turbulent cube has a child whose excess mass produces the contrast
\begin{equation*}
 2^d\mu(Q^\varepsilon)-\mu(Q).
\end{equation*}
The admissible range $\beta<d$ corresponds to the possibility of choosing $\alpha\in(-d,0)$ and $\rho\in(0,1)$ so that $\alpha\rho+\beta<0$. Translated dyadic grids are averaged in order to avoid any assumption concerning mass on the boundaries of dyadic cubes.

The first application is spectral. Let $M(H^1_{\rm at,0}(\T^d))$ denote the unital Banach algebra of Fourier multiplier operators on the complex atomic space. A multiplier $T$ with symbol $m(k)$, $k\in\mathbb Z^d\setminus\{0\}$, is said to have a natural spectrum if
\begin{equation*}
 \sigma_{M(H^1_{\rm at,0}(\T^d))}(T)
 =\overline{\{m(k):k\in\mathbb Z^d\setminus\{0\}\}}.
\end{equation*}
The closure is part of the definition. For a Rajchman measure $\mu$ on $\T^d$ put
\begin{equation}\label{eq:rho-d}
 \rho_{\beta,d}(\mu)
 =\limsup_{m\to\infty}c^{(d)}_\beta(\mu^{*m})^{1/m}
\end{equation}
and
\begin{equation*}
 \|\widehat\mu\|_{\infty,0}
 =\sup_{k\in\mathbb Z^d\setminus\{0\}}|\widehat\mu(k)|.
\end{equation*}
We shall refer to \eqref{eq:rho-d} as the $\beta$-fractal survival radius. When $d=1$ we write $\rho_\beta(\mu)=\rho_{\beta,1}(\mu)$.
The main theorem immediately yields the following multidimensional
criterion for non-natural spectrum.

\begin{cor}[Fractal survival criterion]\label{cor:survival}
Let $d\ge1$ and let $\mu$ be a complex Rajchman measure on $\T^d$. If, for
some $0<\beta<d$,
\begin{equation*}
 \rho_{\beta,d}(\mu)>\|\widehat\mu\|_{\infty,0},
\end{equation*}
then
\begin{equation*}
 \sigma_{M(H^1_{\rm at,0}(\T^d))}(T_\mu)
 \ne
 \overline{\{\widehat\mu(k):k\in\mathbb Z^d\setminus\{0\}\}}.
\end{equation*}
\end{cor}

A particularly transparent consequence is obtained when the whole measure
is carried by a set of lower Minkowski dimension zero.

\begin{cor}[Zero-dimensional Rajchman measures]\label{cor:multidim-zafran}
Let $d\ge1$, and let $\mu$ be a positive Rajchman probability measure on
$\T^d$ carried by a Borel set of lower Minkowski dimension zero. Then
\begin{equation*}
 \sigma_{M(H^1_{\rm at,0}(\T^d))}(T_\mu)
 \ne
 \overline{\{\widehat\mu(k):k\in\mathbb Z^d\setminus\{0\}\}}.
\end{equation*}
\end{cor}

These are genuinely multidimensional statements. Tensor products of
Bluhm's measure give concrete examples for every $d$. When $d=1$, for real measures the
standard identification of the atomic real Hardy space with
$\operatorname{Re}H^1(\T)$ and the passage
$f\mapsto f+i\widetilde f$ also give the corresponding norm estimate and
spectral conclusions on the analytic space $H^1_0(\T)$. We record this
one-dimensional specialization only after the proof of the main theorem.

In dimension one we also obtain a frequency-localized form. For every $\tau>1$ there is an arbitrarily sparse sequence $(k_j)$ such that every frequency set containing all the intervals $[k_j,bk_j^\tau]\cap\mathbb Z$ fails the natural-spectrum property. A single zero-mass Rajchman measure witnesses the failure for all such enlargements, even when the zero frequency is added.

Section \ref{sec:walsh} develops a second group of consequences which is one-dimensional and logically independent of the multiplier estimate. Applying the binary form of Lemma \ref{mrl} directly to a positive measure carried by a set of lower Minkowski dimension less than $\beta<1/2$, we show that, for every $\lambda>\beta/(1-\beta)$, the sum of the $\lfloor 2^{\lambda n}\rfloor$ largest Haar--Stieltjes coefficients at level $n$ is bounded from below along infinitely many levels. Thus a measure concentrated on a sufficiently thin carrier cannot have all of its Haar blocks uniformly compressible.

An elementary Walsh--Haar change of basis, contractive on both $\ell^1$ and $\ell^\infty$, transfers the same lower estimate to the corresponding Walsh--Stieltjes blocks. This yields a Walsh uncertainty principle and, conversely, a geometric consequence: if the sums of the $\lfloor2^{\lambda n}\rfloor$ largest coefficients in the successive Walsh blocks tend to zero, then every full carrier of the measure has lower Minkowski dimension at least $\lambda/(1+\lambda)$. These results explain in detail how the mountain-river combinatorics controls not only multiplier norms but also the distribution of dyadic Fourier mass.

\section{Dyadic combinatorics on \texorpdfstring{$\T^d$}{the torus}}\label{sec:dyadic-d}
For $u\in\T^d$ and $n\ge0$, let
\begin{equation*}
 \mathcal D_n^u
 =\bigl\{u+2^{-n}(j+[0,1)^d)\pmod{\mathbb Z^d}:
 j\in\{0,\ldots,2^n-1\}^d\bigr\}.
\end{equation*}
The families $(\mathcal D_n^u)_{n\ge0}$ form a nested translated dyadic grid. If $Q\in\mathcal D_n^u$, its $2^d$ children in $\mathcal D_{n+1}^u$ are denoted by $Q^\varepsilon$, $\varepsilon\in\{0,1\}^d$.

Let $\nu$ be a finite positive measure. We write $m(Q)=\nu(Q)$ and regard $m$ as a non-negative martingale on the $2^d$-ary tree. For fixed $k$, let $M_d(\beta,k;u)$ be the class of such measures for which
\begin{equation*}
 \#\{Q\in\mathcal D_k^u:m(Q)>0\}\le2^{\beta k}.
\end{equation*}
Fix $-d<\alpha<0$. A cube $Q$ is called turbulent if
\begin{equation}\label{eq:turbulent-d}
 m(Q^\varepsilon)<2^\alpha m(Q)
 \qquad\text{for every }\varepsilon\in\{0,1\}^d.
\end{equation}
Let $\Xi_n^u$ denote the turbulent cubes in $\mathcal D_n^u$.

\begin{lem}[$2^d$-ary mountain river lemma]\label{mrl}\label{lem:mrl-d}
Let $0<\beta<d$, $-d<\alpha<0$, and $0<\rho<1$ satisfy
\begin{equation*}
 \alpha\rho+\beta<0.
\end{equation*}
There are $\theta\in(0,1)$ and $k_0$ such that for every $u\in\T^d$, every $k\ge k_0$, and every $\nu\in M_d(\beta,k;u)$,
\begin{equation}\label{eq:mrl-average}
 \frac1{|I_k|}\sum_{n\in I_k}\sum_{Q\in\Xi_n^u}m(Q)
 <\rho\|\nu\|,
 \qquad I_k=\{n\in\mathbb Z_{\ge0}:\theta k<n<k\}.
\end{equation}
In particular, at some $n\in I_k$ the turbulent cubes carry less than $\rho\|\nu\|$.
\end{lem}

\begin{proof}
Normalize $\|\nu\|=1$. Choose $\rho'<\rho$ so close to $\rho$ that $\alpha\rho'+\beta<0$. Next choose
\begin{equation*}
 0<\theta<1-\frac{\beta}{|\alpha|\rho'}.
\end{equation*}
Put $r=\lfloor\theta k\rfloor$. For a terminal cube $R\in\mathcal D_k^u$, let $S_r(R)$ be the number of turbulent vertices on the branch joining level $r+1$ to level $k-1$. Double counting gives
\begin{equation}\label{eq:branch-count-d}
 \sum_{n=r+1}^{k-1}\sum_{Q\in\Xi_n^u}m(Q)
 =\sum_{R\in\mathcal D_k^u}m(R)S_r(R).
\end{equation}
Split the right-hand side according as
\begin{equation*}
 S_r(R)\le\rho'(k-r-1)
 \quad\text{or}\quad
 S_r(R)>\rho'(k-r-1).
\end{equation*}
The first part is at most $\rho'(k-r-1)$. If $R$ belongs to the second part, then repeated use of \eqref{eq:turbulent-d} along the branch gives
\begin{equation*}
 m(R)\le2^{\alpha\rho'(k-r-1)}.
\end{equation*}
There are at most $2^{\beta k}$ terminal cubes with positive mass. Hence the second part is bounded by
\begin{equation*}
 k\,2^{\beta k+\alpha\rho'(k-r-1)}.
\end{equation*}
Since $r/k\to\theta$ and
\begin{equation*}
 \beta+\alpha\rho'(1-\theta)<0,
\end{equation*}
this term is $o(k-r-1)$. Thus, for all sufficiently large $k$, division of \eqref{eq:branch-count-d} by $k-r-1$ gives
\begin{equation*}
 \frac1{k-r-1}
 \sum_{n=r+1}^{k-1}\sum_{Q\in\Xi_n^u}m(Q)<\rho.
\end{equation*}
This is \eqref{eq:mrl-average}.
\end{proof}

The following elementary observation permits the use of arbitrary translated grids. The small covering number is preserved.

\begin{lem}[Thin coverings in translated grids]\label{lem:translated-cover}
Let $A\subset\T^d$ and $\ldimM A<\beta$. Then there are arbitrarily large integers $k$ such that, for every $u\in\T^d$, the set $A$ meets fewer than $2^{\beta k}$ cubes of $\mathcal D_k^u$.
\end{lem}

\begin{proof}
Choose $\beta'$ with $\ldimM A<\beta'<\beta$. There are arbitrarily small numbers $\delta$ for which
\begin{equation*}
 N(A,\delta)\le \delta^{-\beta'}.
\end{equation*}
For each such $\delta$, choose $k$ so that
\begin{equation*}
 2^{-k-1}<\delta\le2^{-k}.
\end{equation*}
A cube of side $\delta$ meets at most $C_d$ cubes of any translated grid $\mathcal D_k^u$. Hence $A$ meets at most
\begin{equation*}
 C_d\delta^{-\beta'}\le C_d2^{\beta'(k+1)}
\end{equation*}
level-$k$ cubes. Since $\beta'<\beta$, this is smaller than $2^{\beta k}$ for all sufficiently large members of the resulting sequence of integers $k$.
\end{proof}

For $0<q<1/8$, let $B_n^u(q)$ be the set of points whose $\ell^\infty$-distance from the boundary of a cube in $\mathcal D_{n+1}^u$ is less than $q2^{-n}$.

\begin{lem}[Averaging away the dyadic boundaries]\label{lem:boundary-average}
For every finite positive measure $\nu$ on $\T^d$,
\begin{equation}\label{eq:boundary-average}
 \int_{\T^d}\nu(B_n^u(q))\,du\le4dq\|\nu\|.
\end{equation}
\end{lem}

\begin{proof}
Fix $x\in\T^d$. For one coordinate, the set of shifts $u$ for which $x$ lies within $q2^{-n}$ of a level-$(n+1)$ grid hyperplane has normalized measure at most $4q$. Taking the union over the $d$ coordinates gives at most $4dq$. Integrating first in $u$ and then in $x$ proves \eqref{eq:boundary-average}.
\end{proof}

\begin{lem}[Geometric testing estimate]\label{lem:geometric-testing}
Fix $u\in\T^d$, $n\ge0$, $0<q<1/8$, and $\varepsilon\in\{0,1\}^d$. Let $\mathcal G\subset\mathcal D_n^u$, put $L=2^{-n}$ and $b=2^d$, and define
\begin{equation*}
 D_\varepsilon=\frac L2\bigl(\varepsilon+[0,1)^d\bigr),
 \qquad C=[-L,0)^d,
 \qquad C_\varepsilon=-D_\varepsilon,
\end{equation*}
\begin{equation*}
 h(x)=L^{-d}\bigl(b\chi_{C_\varepsilon}(x)-\chi_C(x)\bigr).
\end{equation*}
The choice of half-open representatives on the boundaries is immaterial for the integrals below. For each $Q\in\mathcal G$, choose a lift $x_Q+[0,L)^d$ of $Q$ to $\mathbb R^d$ and set
\begin{equation*}
 E_Q=x_Q+[0,qL)^d,
 \qquad E=\bigcup_{Q\in\mathcal G}E_Q
\end{equation*}
after projection to $\T^d$. Then the sets $E_Q$ are disjoint and, for every finite positive measure $\nu$,
\begin{equation}\label{eq:geometric-testing}
 \int_E h*\nu(t)\,dt
 \ge q^d\sum_{Q\in\mathcal G}
 \bigl(b\nu(Q^\varepsilon)-\nu(Q)\bigr)
 -C_dq^d\nu(B_n^u(q)).
\end{equation}
Moreover, if every $Q\in\mathcal G$ meets a set $A\subset\T^d$, then for every finite positive measure $\eta$,
\begin{equation}\label{eq:geometric-remainder}
 \int_E|h*\eta(t)|\,dt
 \le 2\eta\bigl(\{x:\operatorname{dist}_\infty(x,A)<3L\}\bigr).
\end{equation}
All distances are torus distances, so the assertions are unaffected by cubes crossing a fundamental-domain boundary.
\end{lem}

\begin{proof}
The sets $E_Q$ lie inside the pairwise disjoint cubes $Q$. Work first with the chosen lifts and then project periodically. For $x\notin B_n^u(q)$ and almost every $t\in E_Q$, translation by the vector $t-x_Q\in[0,qL)^d$ cannot move $x$ across a boundary of a level-$(n+1)$ cube. Consequently
\begin{equation*}
 h(t-x)=
 \begin{cases}
 (b-1)L^{-d},&x\in Q^\varepsilon,\\
 -L^{-d},&x\in Q\setminus Q^\varepsilon,\\
 0,&x\notin Q.
 \end{cases}
\end{equation*}
Since $|E_Q|=q^dL^d$, integration in $t$ gives
\begin{equation*}
 \sum_{Q\in\mathcal G}\int_{E_Q}h(t-x)\,dt
 =q^d\sum_{Q\in\mathcal G}
 \bigl(b\chi_{Q^\varepsilon}(x)-\chi_Q(x)\bigr)
\end{equation*}
whenever $x\notin B_n^u(q)$.

For arbitrary $x$, a term indexed by $Q$ can be non-zero only when
\begin{equation*}
 x\in E_Q-C=x_Q+[0,(1+q)L)^d.
\end{equation*}
These enlarged cubes have overlap bounded by a constant depending only on $d$, also after periodic projection. As $\|h\|_\infty\le(b-1)L^{-d}$, both sides of the preceding pointwise identity are bounded in absolute value by $C_dq^d$. Their difference is therefore bounded by $C_dq^d\chi_{B_n^u(q)}(x)$. Integration with respect to $\nu$ and Fubini's theorem prove \eqref{eq:geometric-testing}.

Finally, if $t\in E_Q$ and $h(t-x)\ne0$, then $x\in E_Q-C$. If $Q$ meets $A$, the torus $\ell^\infty$-distance from $x$ to $A$ is less than $3L$. Hence, by Fubini and $\|h\|_1=2(b-1)/b<2$,
\begin{equation*}
 \int_E|h*\eta(t)|\,dt
 \le\int\!\int_E|h(t-x)|\,dt\,d\eta(x)
 \le2\eta\bigl(\{x:\operatorname{dist}_\infty(x,A)<3L\}\bigr),
\end{equation*}
which is \eqref{eq:geometric-remainder}.
\end{proof}

\section{The multidimensional Hardy estimate}\label{sec:main-proof}
We prove a scale-by-scale statement which immediately implies the main theorem.

\begin{prop}[Scale-by-scale form on $\T^d$]\label{prop:scale-d}
Let $d\ge1$ and $0<\beta<d$. There are constants $c_0(d,\beta)>0$ and $A_d>0$ such that for every real measure $\mu\in M(\T^d)$ with $c^{(d)}_\beta(\mu)>0$, there are infinitely many integers $n$ and real test functions $h_n\in H^1_{\rm at,0}(\T^d)$ satisfying
\begin{equation}\label{eq:scale-d}
 \|h_n\|_{H^1_{\rm at,0}}\le A_d,
 \qquad
 \|h_n*\mu\|_1\ge c_0(d,\beta)c^{(d)}_\beta(\mu).
\end{equation}
Each $h_n$ is supported on a cube of side $2^{-n}$ and is a constant multiple, bounded in terms of $d$ only, of the difference between the indicator of one child and the average over its parent.
\end{prop}

\begin{proof}
Put $c=c^{(d)}_\beta(\mu)$. Choose
\begin{equation}\label{eq:parameter-choice-d}
 \beta<a<d,
 \qquad \frac\beta a<\rho<1,
 \qquad \alpha=-a.
\end{equation}
Then $-d<\alpha<0$, $\alpha\rho+\beta<0$, and
\begin{equation}\label{eq:positive-contrast-d}
 2^{d+\alpha}-1=2^{d-a}-1>0.
\end{equation}
Let $\theta$ be supplied by Lemma \ref{mrl}. Fix a small $\tau>0$, eventually chosen as a sufficiently small multiple of $c$.

By the definition of $c$, choose a Borel set $B$ with $\ldimM B<\beta$ and $|\mu|(B)>c-\tau$. Let $\mu=\mu^+-\mu^-$ be the Jordan decomposition, and choose a Borel Hahn decomposition $\T^d=P\cup N$ with $\mu^-|_P=0$ and $\mu^+|_N=0$. One of $\mu^+(B)$ and $\mu^-(B)$ is at least $(c-\tau)/2$. Replacing $\mu$ by $-\mu$ if necessary, we may suppose that $\mu^+(B)\ge(c-\tau)/2$. By regularity choose a compact set $A\subset B\cap P$ such that
\begin{equation}\label{eq:one-sign-d}
 \mu|_A\ge0,
 \qquad \ldimM A<\beta,
 \qquad \mu(A)>\frac c2-\tau.
\end{equation}
Put
\begin{equation*}
 \nu=\mu|_A,
 \qquad \sigma=\mu-\nu,
 \qquad \eta=|\sigma|=|\mu||_{\T^d\setminus A}.
\end{equation*}
Then $\nu$ is positive, $\|\nu\|>c/2-\tau$, and $\eta(A)=0$.

By Lemma \ref{lem:translated-cover}, there are arbitrarily large $k$ such that $A$ meets fewer than $2^{\beta k}$ cubes of every translated grid $\mathcal D_k^u$. We take such a $k$, sufficiently large for Lemma \ref{mrl}, and also so large that
\begin{equation}\label{eq:remainder-neighborhood-d}
 \eta\bigl(\{x:\operatorname{dist}_\infty(x,A)<C_d2^{-\theta k}\}\bigr)<\tau.
\end{equation}
This is possible because $A$ is compact and $\eta(A)=0$.

Put $b=2^d$ and
\begin{equation*}
 \kappa=\frac{b2^\alpha-1}{b}>0.
\end{equation*}
For fixed $u$ and $n\in I_k$, let $\mathcal N_n^u$ be the non-turbulent cubes of positive $\nu$-mass. For every $Q\in\mathcal N_n^u$, choose a child $Q^{\varepsilon(Q)}$ such that
\begin{equation*}
 \nu(Q^{\varepsilon(Q)})\ge2^\alpha\nu(Q).
\end{equation*}
For $\varepsilon\in\{0,1\}^d$, let $G_{n,u}^\varepsilon$ be the cubes for which $\varepsilon(Q)=\varepsilon$. By \eqref{eq:positive-contrast-d},
\begin{equation}\label{wzn}
 \max_\varepsilon
 \sum_{Q\in G_{n,u}^\varepsilon}
 \bigl(b\nu(Q^\varepsilon)-\nu(Q)\bigr)
 \ge \kappa\sum_{Q\in\mathcal N_n^u}\nu(Q).
\end{equation}

Lemma \ref{mrl} implies that the average over $n\in I_k$ of $\sum_{Q\in\mathcal N_n^u}\nu(Q)$ is larger than $(1-\rho)\|\nu\|$, uniformly in $u$. Combining this with \eqref{wzn} and Lemma \ref{lem:boundary-average} gives the explicit averaged estimate
\begin{equation*}
\begin{aligned}
 &\mathbb E_{u,n}\!\left[
 \max_{\varepsilon}\sum_{Q\in G_{n,u}^{\varepsilon}}
 \bigl(b\nu(Q^\varepsilon)-\nu(Q)\bigr)
 -C_d\nu(B_n^u(q))
 \right]\\
 &\hspace{35mm}\ge \kappa(1-\rho)\|\nu\|-4C_ddq\|\nu\|.
\end{aligned}
\end{equation*}
where $u$ is averaged with respect to normalized Haar measure and $n$ uniformly over $I_k$. Choose $q=q(d,\alpha,\rho)>0$ so small that the right-hand side is at least $\frac12\kappa(1-\rho)\|\nu\|$. We can then select a common triple $u$, $n\in I_k$, and $\varepsilon$ such that

\begin{equation}\label{eq:good-grid-level-d}
 \sum_{Q\in G_{n,u}^\varepsilon}
 \bigl(b\nu(Q^\varepsilon)-\nu(Q)\bigr)
 -C_d\nu(B_n^u(q))
 \ge \frac{\kappa(1-\rho)}2\|\nu\|.
\end{equation}

Apply Lemma \ref{lem:geometric-testing} with $\mathcal G=G_{n,u}^\varepsilon$. Thus, with $L=2^{-n}$,
\begin{equation}\label{eq:test-atom-d}
 h_n(x)=L^{-d}\bigl(b\chi_{C_\varepsilon}(x)-\chi_C(x)\bigr),
\end{equation}
and the testing set $E$ is the union of the cubes $E_Q$ from that lemma. The function $h_n$ has integral zero, is supported on a cube of side $L$, and satisfies
\begin{equation*}
 \|h_n\|_\infty\le(b-1)L^{-d}.
\end{equation*}
Thus $(b-1)^{-1}h_n$ is an atom and
\begin{equation}\label{eq:atom-norm-d}
 \|h_n\|_{H^1_{\rm at,0}}\le b-1.
\end{equation}
By \eqref{eq:geometric-testing} and \eqref{eq:good-grid-level-d},
\begin{align}
 \int_E h_n*\nu(t)\,dt
 &\ge q^d\sum_{Q\in G_{n,u}^\varepsilon}
 \bigl(b\nu(Q^\varepsilon)-\nu(Q)\bigr)
 -C_dq^d\nu(B_n^u(q))\notag\\
 &\ge c_1(d,\beta)\|\nu\|.\label{eq:conv-lower-d}
\end{align}
In particular, $\|h_n*\nu\|_1\ge c_1(d,\beta)\|\nu\|$.

Every cube $Q$ in the construction has positive $\nu$-mass and hence meets $A$. Since $n>\theta k$, we have $3L\le C_d2^{-\theta k}$. Applying \eqref{eq:geometric-remainder} to the positive measure $\eta$ and then \eqref{eq:remainder-neighborhood-d}, we obtain
\begin{equation*}
 \int_E |h_n*\sigma(t)|\,dt
 \le \int_E (|h_n|*\eta)(t)\,dt
 \le2\tau.
\end{equation*}
Consequently,
\begin{equation*}
 \|h_n*\mu\|_1
 \ge c_1(d,\beta)\left(\frac c2-\tau\right)-C_d\tau.
\end{equation*}
Choosing $\tau$ to be a sufficiently small multiple of $c$ proves \eqref{eq:scale-d}. Since the admissible covering scales $k$ are arbitrarily large and $n>\theta k$, the testing scales $n$ form an infinite set.
\end{proof}

\begin{proof}[Proof of Theorem \ref{dwa}]
First suppose that $\mu$ is real. If $c^{(d)}_\beta(\mu)=0$, there is nothing to prove; otherwise Proposition \ref{prop:scale-d} and \eqref{eq:atom-norm-d} give \eqref{eq:main-multidim}.

For a general complex measure put $c=c^{(d)}_\beta(\mu)$. If $c=0$, there is nothing to prove. Write $d\mu=\phi\,d|\mu|$, where $|\phi|=1$ almost everywhere, and choose a Borel set $B$ with $\ldimM B<\beta$ and $|\mu|(B)>c/2$. Partition the unit circle into eight half-open arcs of equal length. One phase sector carries at least $|\mu|(B)/8>c/16$. If $\vartheta$ is the midpoint of that sector, let $B_0\subset B$ be its inverse image under $\phi$ and put
\begin{equation*}
 \lambda=\operatorname{Re}(e^{-i\vartheta}\mu).
\end{equation*}
On $B_0$ we have $d\lambda\ge \cos(\pi/8)\,d|\mu|$, and therefore
\begin{equation*}
 c^{(d)}_\beta(\lambda)
 \ge |\lambda|(B_0)
 \ge \lambda(B_0)
 \ge \frac{\cos(\pi/8)}{16}\,c.
\end{equation*}
Apply the real case to $\lambda$. The resulting test functions $h_n$ are real, occur at infinitely many scales, and satisfy
\begin{equation*}
 \|h_n*\mu\|_1
 \ge \|\operatorname{Re}(e^{-i\vartheta}h_n*\mu)\|_1
 =\|h_n*\lambda\|_1.
\end{equation*}
This proves the complex case, with a changed constant.

Finally, translations are isometries of the atomic space and
\begin{equation*}
 T_\mu f=\int_{\T^d}\tau_y f\,d\mu(y)
\end{equation*}
in the Bochner sense. Hence $\|T_\mu f\|_{H^1}\le\|\mu\|\,\|f\|_{H^1}$. Since the Hardy norm dominates the $L^1$ norm, the same lower estimate holds for the operator norm on $H^1_{\rm at,0}(\T^d)$.
\end{proof}

\begin{rem}[The one-dimensional analytic specialization]\label{rem:analytic-specialization}
For $d=1$, the atomic space $H^1_{\rm at,0}(\T)$ is equivalent to the
usual real Hardy space $\operatorname{Re}H^1(\T)$. If $h$ is one of the
real test functions from Proposition \ref{prop:scale-d}, then

If $F=h+i\widetilde h$, then $F\in H^1_0(\T)$ and there is an absolute constant $C_H\ge1$, independent of $h$, such that
\begin{equation*}
 C_H^{-1}\|h\|_{\operatorname{Re}H^1}
 \le \|F\|_{H^1}
 \le C_H\|h\|_{\operatorname{Re}H^1}.
\end{equation*}
Moreover, for every real measure $\mu$ one has
\begin{equation*}
 \sup_{n\ge1}|\widehat\mu(n)|
 =\sup_{n\ne0}|\widehat\mu(n)|,
\end{equation*}
because $\widehat\mu(-n)=\overline{\widehat\mu(n)}$.

For a real measure $\mu$,
\begin{equation*}
 \operatorname{Re}(T_\mu F)=T_\mu h,
\end{equation*}
so the $d=1$ case of Theorem \ref{dwa} also gives
\begin{equation*}
 \|T_\mu:H^1_0(\T)\to H^1_0(\T)\|
 \ge C(\beta)c_\beta(\mu),
 \qquad 0<\beta<1.
\end{equation*}
Consequently, for real Rajchman measures the conclusion of Corollary \ref{cor:survival}, and for positive measures the conclusion of Corollary \ref{cor:multidim-zafran}, have their analytic counterparts on $H^1_0(\T)$. This is a one-dimensional by-product of the
multidimensional real-variable theorem.
\end{rem}

For the frequency-localized one-dimensional application below we record the additional regularity of the testing atoms supplied by the proof.

\begin{prop}[One-dimensional scale-by-scale form]\label{prop:scale}
For every $0<\beta<1$ there are constants $c_0(\beta)>0$ and $A>0$ such that for every real $\mu\in M(\T)$ with $c_\beta(\mu)>0$, there are infinitely many $n$ and real functions $h_n$ satisfying
\begin{equation*}
 \|h_n\|_{\mathrm{Re}H^1}\le A,
 \qquad
 \|h_n*\mu\|_1\ge c_0(\beta)c_\beta(\mu),
\end{equation*}
and
\begin{equation*}
 \operatorname{supp}h_n\subset I_n,
 \quad |I_n|=2^{-n},
 \quad \int h_n=0,
\end{equation*}
\begin{equation*}
 \|h_n\|_1\le A,
 \qquad \|h_n\|_\infty\le A2^n,
 \qquad \operatorname{Var}(h_n)\le A2^n.
\end{equation*}
\end{prop}

\begin{proof}
Take $d=1$ in Proposition \ref{prop:scale-d}. The test function \eqref{eq:test-atom-d} is, up to translation and sign, the normalized Haar function on an interval of length $2^{-n}$. Hence its real Hardy norm and $L^1$ norm are uniformly bounded, its $L^\infty$ norm is $O(2^n)$, and its total variation is $O(2^n)$.
\end{proof}

\section{Spectral applications}\label{sec:applications}
We now prove the multidimensional spectral corollaries stated in the introduction. Spectra are taken in the multiplier algebra $M(H^1_{\rm at,0}(\T^d))$, and natural spectrum always has the meaning fixed in the introduction, including closure of the symbol range. A measure $\mu$ on $\T^d$ is Rajchman if
$\widehat\mu(k)\to0$ as $|k|\to\infty$ in $\mathbb Z^d$.

The following elementary estimates will be used for convolution powers.

\begin{lem}[Closure invariance of lower Minkowski dimension]\label{lem:closure-minkowski}
For every $E\subset\T^d$,
\begin{equation*}
 \ldimM\overline E=\ldimM E.
\end{equation*}
\end{lem}
\begin{proof}
The inequality $\ldimM E\le\ldimM\overline E$ is immediate. Conversely, if $E$ is covered by $N(E,\delta)$ half-open cubes of side $\delta$, then $\overline E$ is contained in the union of their closures. Each such closed cube can be covered by at most $C_d$ half-open cubes of side $2\delta$, where $C_d$ depends only on the dimension and on the fixed torus convention. Hence
\begin{equation*}
 N(\overline E,2\delta)\le C_dN(E,\delta).
\end{equation*}
Taking logarithms and the lower limit proves the reverse inequality.
\end{proof}

\begin{lem}\label{lem:sumsets}
Let $A\subset\T^d$ and $m\in\N$. Then
\begin{equation*}
 \ldimM(mA)\le m\,\ldimM A,
\end{equation*}
Here $mA=A+\cdots+A$.\par In particular, $\ldimM A=0$ implies $\ldimM(mA)=0$ for every $m$.
\end{lem}

\begin{proof}
If $A$ is covered by $N(A,\delta)$ cubes of side $\delta$, then $mA$ is covered by at most $N(A,\delta)^m$ cubes of side $m\delta$. Hence
\begin{equation*}
 N(mA,m\delta)\le N(A,\delta)^m.
\end{equation*}
Taking logarithms and the lower limit along a sequence realizing $\ldimM A$ proves the assertion. See also \citep[Proposition 1, p.~71]{BS}.
\end{proof}

\begin{proof}[Proof of Corollary \ref{cor:survival}]
By Theorem \ref{dwa},
\begin{equation*}
 \|T_\mu^m\|
 =\|T_{\mu^{*m}}\|
 \ge C(d,\beta)c^{(d)}_\beta(\mu^{*m}).
\end{equation*}
Taking $m$-th roots and the upper limit gives
\begin{equation*}
 r(T_\mu)\ge\rho_{\beta,d}(\mu).
\end{equation*}
If the spectrum were natural, its spectral radius would equal
$\|\widehat\mu\|_{\infty,0}$, a contradiction.
\end{proof}

\begin{proof}[Proof of Corollary \ref{cor:multidim-zafran}]

Let $A$ be a Borel carrier of $\mu$ with $\ldimM A=0$ and put $K=\overline A$. By Lemma \ref{lem:closure-minkowski}, $K$ is a compact carrier with $\ldimM K=0$. Hence $mK$ is compact, and therefore Borel, while Lemma \ref{lem:sumsets} gives $\ldimM(mK)=0$. Since $\mu^{*m}$ is carried by $mK$ and is a probability measure,
\begin{equation*}
 c^{(d)}_\beta(\mu^{*m})=1
 \qquad(m\ge1).
\end{equation*}
Thus $\rho_{\beta,d}(\mu)=1$ for every $0<\beta<d$.

Moreover,
\begin{equation*}
 \|\widehat\mu\|_{\infty,0}<1.
\end{equation*}
Indeed, equality at a non-zero frequency $k_0$ would force the
corresponding character to be constant $\mu$-almost everywhere and hence
$|\widehat\mu(mk_0)|=1$ for all $m$, contradicting the Rajchman
property. Thus every individual non-zero coefficient has modulus strictly below $1$. Since the coefficients tend to zero at infinity, their supremum is the maximum of finitely many of them together with a tail of arbitrarily small modulus; hence the supremum is also strictly below $1$. Corollary \ref{cor:survival} applies.
\end{proof}

\begin{lem}[The diffuse measures form a convolution ideal]\label{lem:diffuse-ideal}
Let $0<\beta<d$ and let $\eta\in M(\T^d)$ satisfy $c^{(d)}_\beta(\eta)=0$. Then, for every $\sigma\in M(\T^d)$,
\begin{equation*}
 c^{(d)}_\beta(\eta*\sigma)=0.
\end{equation*}
More precisely, if $B\subset\T^d$ is Borel and $\ldimM B<\beta$, then $(|\eta|*|\sigma|)(B)=0$.
\end{lem}

\begin{proof}
Every translate $B-x$ has the same lower Minkowski dimension as $B$. Therefore $|\eta|(B-x)=0$ for all $x$, and
\begin{equation*}
 (|\eta|*|\sigma|)(B)
 =\int_{\T^d}|\eta|(B-x)\,d|\sigma|(x)=0.
\end{equation*}
Since $|\eta*\sigma|\le|\eta|*|\sigma|$, the conclusion follows.
\end{proof}

\begin{cor}[Thin component plus diffuse remainder]\label{cor:noise}
Let $0<\beta<d$. Let $\nu$ be a real Rajchman measure on $\T^d$ carried by a Borel set $A$ such that
\begin{equation*}
 \ldimM(mA)<\beta\qquad(m\ge1),
\end{equation*}
and let $\eta$ be a real Rajchman measure satisfying $c^{(d)}_\beta(\eta)=0$. Put
\begin{equation*}
 \mu=a\nu+\eta,
 \qquad a\in\mathbb R,
\end{equation*}
and let
\begin{equation*}
 r_M(\nu)=\lim_{m\to\infty}\|\nu^{*m}\|^{1/m}
\end{equation*}
be the spectral radius of $\nu$ in $M(\T^d)$. If
\begin{equation*}
 |a|r_M(\nu)>\|\widehat\mu\|_{\infty,0},
\end{equation*}
then $T_\mu$ has non-natural spectrum on $H^1_{\rm at,0}(\T^d)$. A sufficient explicit condition is
\begin{equation*}
 |a|\|\widehat\nu\|_{\infty,0}
 +\|\widehat\eta\|_{\infty,0}
 <|a|r_M(\nu).
\end{equation*}
If $\nu$ is a positive probability measure, then $r_M(\nu)=1$.
\end{cor}

\begin{proof}

Put $K=\overline A$. Since $A^m$ is dense in $K^m$ and addition is continuous, compactness gives
\begin{equation*}
 mK=\overline{mA}.
\end{equation*}
Lemmas \ref{lem:closure-minkowski} and \ref{lem:sumsets}, together with the hypothesis on $mA$, imply that $mK$ is a compact Borel set with $\ldimM(mK)<\beta$. Expand $\mu^{*m}$. The term containing no factor $\eta$ is $a^m\nu^{*m}$; denote the sum of the other terms by $R_m$. Every summand of $R_m$ contains a factor $\eta$, so Lemma \ref{lem:diffuse-ideal} gives $c^{(d)}_\beta(R_m)=0$ and hence $|R_m|(mK)=0$. Since $\nu^{*m}$ is carried by $mK$, the restrictions of the measures satisfy
\begin{equation*}
 \mu^{*m}|_{mK}=a^m\nu^{*m}|_{mK},
\end{equation*}
and therefore
\begin{equation*}
 c^{(d)}_\beta(\mu^{*m})
 \ge|\mu^{*m}|(mK)
 =|a|^m\|\nu^{*m}\|.
\end{equation*}

Hence $\rho_{\beta,d}(\mu)\ge|a|r_M(\nu)$, and Corollary \ref{cor:survival} applies. The explicit sufficient condition follows from
\begin{equation*}
 \|\widehat\mu\|_{\infty,0}
 \le |a|\|\widehat\nu\|_{\infty,0}
 +\|\widehat\eta\|_{\infty,0}.
\end{equation*}
\end{proof}

\begin{rem}
The condition $c^{(d)}_\beta(\eta)=0$ says that, once the thin component has been selected, the remainder has no further mass on sets of lower Minkowski dimension below $\beta$. Any thin part still present in the remainder should simply be absorbed into $\nu$. Absolute continuity is only the simplest example of this diffuse property.
\end{rem}

We recall why the above statements are non-vacuous. Bluhm \citep{B} constructed a positive Rajchman measure carried by a suitable Cantor subset of the Liouville numbers. Here $\mathbb P_M$ denotes the set of prime numbers $p$ with $M\le p<2M$, and $\|y\|=\min_{j\in\mathbb Z}|y-j|$ denotes the distance from $y$ to the nearest integer. In the notation of that construction, choose integers
\begin{equation*}
 M_1<2M_1<M_2<2M_2<\cdots
\end{equation*}
and put
\begin{equation*}
 S_\infty
 =\bigcap_{k=1}^\infty
 \bigcup_{p\in\mathbb P_{M_k}}
 \{x\in[0,1]:\|px\|\le p^{-1-k}\}.
\end{equation*}
For a suitable choice of $(M_k)$, the set $S_\infty$ supports a positive Rajchman probability measure $\nu$. For fixed $k$ and $p\in\mathbb P_{M_k}$, the corresponding set is a union of at most $p$ intervals of length at most $2p^{-2-k}$. Hence $S_\infty$ has a cover by at most $CM_k^2$ intervals of length at most $CM_k^{-k-2}$, and therefore
\begin{equation*}
 \frac{\log N(S_\infty,CM_k^{-k-2})}
 {\log(M_k^{k+2}/C)}
 \le
 \frac{2\log M_k+O(1)}{(k+2)\log M_k+O(1)}\longrightarrow0.
\end{equation*}
Thus $\ldimM S_\infty=0$.

\begin{rem}[Tensor products of Bluhm measures]\label{rem:tensor-bluhm}
Let $\nu$ be the measure above and put
\begin{equation*}
 \mu_d=\nu^{\otimes d}.
\end{equation*}
Then
\begin{equation*}
 \widehat{\mu_d}(k_1,\ldots,k_d)
 =\prod_{j=1}^d\widehat\nu(k_j),
\end{equation*}
so $\mu_d$ is Rajchman. It is carried by $S_\infty^d$, and
\begin{equation*}
 N(S_\infty^d,\delta)\le N(S_\infty,\delta)^d,
\end{equation*}
whence $\ldimM(S_\infty^d)=0$. Thus Corollary \ref{cor:multidim-zafran} is non-vacuous in every dimension. The tensor product is essential: the diagonal image of $\nu$ would have Fourier transform equal to one on non-trivial frequency hyperplanes and would not be Rajchman on $\T^d$.
\end{rem}

We can proceed even further. For $\Lambda\subset\mathbb{Z}$, let $L^{1}_{\Lambda}(\T)$ denote the subspace of summable functions whose Fourier coefficients vanish outside $\Lambda$, and let $M(L^{1}_{\Lambda}(\T))$ be its multiplier algebra. We write
\begin{equation*}
 M_0(L^1_\Lambda(\T))
 =\bigl\{(\widehat\mu(n))_{n\in\Lambda}:\mu\in M_0(\T)\bigr\}
 \subset M(L^1_\Lambda(\T))
\end{equation*}
for the image of the Rajchman measures under restriction of their Fourier--Stieltjes transforms to $\Lambda$. A multiplier $S\in M(L^1_\Lambda(\T))$ with symbol $m(n)$, $n\in\Lambda$, is said to have a natural spectrum when $\sigma(S)=\overline{\{m(n):n\in\Lambda\}}$. A refinement of the previous argument allows us to construct a much smaller set $\Lambda$ for which there is a multiplier without a natural spectrum.
\begin{lem}[Frequency localization of the testing atoms]\label{lem:frequency-localization}
Fix $\tau>1$ and choose
\begin{equation*}
        \frac1\tau<\gamma<1.
\end{equation*}
Let $(h_n)$ be the real test functions supplied by Proposition \ref{prop:scale}. There are analytic trigonometric polynomials $p_n$ and positive constants $C_{\tau,\gamma}$ and $b_{\tau,\gamma}$ such that
\begin{equation*}
        \operatorname{spec}p_n\subset [N_n,b_{\tau,\gamma}N_n^\tau],
        \qquad N_n\asymp 2^{\gamma n},
\end{equation*}
\begin{equation*}
        \|p_n\|_1\le C_{\tau,\gamma},
\end{equation*}
and
\begin{equation*}
        \|\Re p_n-h_n\|_1
        \le C_{\tau,\gamma}\bigl(2^{(\gamma-1)n}+2^{(1-\tau\gamma)n}\bigr).
\end{equation*}
In particular, $\|\Re p_n-h_n\|_1\to0$ as $n\to\infty$.
\end{lem}
\begin{proof}
Let $V_m$ denote the de la Vall\'ee Poussin operator, normalized so that its Fourier multiplier is equal to $1$ on $[-m,m]$ and to $0$ outside $[-2m,2m]$.

Its kernels $K_m$ have uniformly bounded $L^1$ norms. Since translations act isometrically on $\operatorname{Re}H^1$, the convolution representation $V_mf=K_m*f$ gives
\begin{equation*}
 \|V_mf\|_{\operatorname{Re}H^1}
 \le \|K_m\|_1\|f\|_{\operatorname{Re}H^1},
\end{equation*}
so the operators are uniformly bounded on $\operatorname{Re}H^1$. We shall also use
\begin{equation*}
 \|K_m'\|_1\le Cm\|K_m\|_1\le Cm,
\end{equation*}
which is the Bernstein inequality for trigonometric polynomials. Finally, the classical $L^1$ Jackson inequality for functions of bounded variation \citep[Vol.~I, p.~114]{Zygmund} gives
\begin{equation*}
 E_m(f)_1:=\inf_{\deg q\le m}\|f-q\|_1
 \le C\frac{\operatorname{Var}(f)}m.
\end{equation*}
Because $V_mq=q$ for $\deg q\le m$ and $V_m$ is uniformly bounded on $L^1$, it follows that
\begin{equation}\label{eq:vp-bv-approx}
 \|f-V_mf\|_1
 \le(1+\|V_m\|_{L^1\to L^1})E_m(f)_1
 \le C\frac{\operatorname{Var}(f)}m.
\end{equation}

Put
\begin{equation*}
        A_n=\lfloor 2^{\gamma n}\rfloor,
        \qquad
        B_n=\lfloor 2^{\tau\gamma n}\rfloor,
\end{equation*}
and define
\begin{equation*}
        g_n=V_{B_n}h_n-V_{A_n}h_n .
\end{equation*}
The function $g_n$ is real and has mean zero. Let $p_n=g_n+i\widetilde g_n$ be the corresponding analytic polynomial. Then $p_n\in H^1_0$, $\Re p_n=g_n$, and its positive spectrum is contained in
\begin{equation*}
        [A_n+1,2B_n]\subset [N_n,b_{\tau,\gamma}N_n^\tau]
\end{equation*}
with $N_n=A_n+1\asymp2^{\gamma n}$. Uniform boundedness on $\operatorname{Re}H^1$ and Proposition \ref{prop:scale} give
\begin{equation*}
        \|p_n\|_1\le C_{\tau,\gamma}.
\end{equation*}

It remains to estimate $h_n-g_n$.  Applying \eqref{eq:vp-bv-approx} to $h_n$ gives
\begin{equation*}
        \|h_n-V_{B_n}h_n\|_1
        \le C\frac{\operatorname{Var}(h_n)}{B_n}
        \le C2^{(1-\tau\gamma)n}.
\end{equation*}
For the low-frequency part, let $x_n$ be any point of $\operatorname{supp}h_n$. Since $\int h_n=0$,
\begin{equation*}
 V_{A_n}h_n(x)
 =\int h_n(t)\bigl(K_{A_n}(x-t)-K_{A_n}(x-x_n)\bigr)\,dt.
\end{equation*}
Therefore, using the translation estimate in $L^1$ and
$|t-x_n|\le 2^{-n}$ on the support of $h_n$,
\begin{align*}
 \|V_{A_n}h_n\|_1
 &\le \int |h_n(t)|\,
       \|\tau_tK_{A_n}-\tau_{x_n}K_{A_n}\|_1\,dt \\
 &\le \|h_n\|_1\,2^{-n}\|K_{A_n}'\|_1
 \le C A_n2^{-n}
 \le C2^{(\gamma-1)n}.
\end{align*}
Combining the two estimates proves the lemma.
\end{proof}

\begin{rem}
The cubic upper endpoint used in the earlier formulation is not intrinsic. The preceding proof gives every power $\tau>1$. The endpoint $\tau=1$ is not reached by this localization argument: removal of the low frequencies requires $A_n2^{-n}\to0$, whereas approximation of the jump discontinuities requires $2^n/B_n\to0$, and hence necessarily $B_n/A_n\to\infty$.
\end{rem}

\begin{prop}\label{prop:localized}
For every $\beta\in(0,1)$ and every $\tau>1$ there are constants
$b=b(\tau)>0$ and $C_{\rm loc}(\beta,\tau)>0$ such that, for every real measure $\mu\in M(\T)$, there are arbitrarily large integers $N$ with the following property: if
\begin{equation*}
        \Lambda=[N,bN^\tau]\cap\mathbb Z,
\end{equation*}
then
\begin{equation*}
\|T_{\mu}:L^{1}_{\Lambda}(\T)\to L^{1}_{\Lambda}(\T)\|
\geq C_{\rm loc}(\beta,\tau)c_{\beta}(\mu).
\end{equation*}
\end{prop}
\begin{proof}
If $c_{\beta}(\mu)=0$, there is nothing to prove. Assume $c_{\beta}(\mu)>0$. By Proposition \ref{prop:scale}, there are infinitely many scales $n$ and test functions $h_n$ for which
\begin{equation*}
        \|T_\mu h_n\|_1\ge c_0(\beta)c_{\beta}(\mu).
\end{equation*}
Choose $\gamma\in(1/\tau,1)$ and let $p_n$ be supplied by Lemma \ref{lem:frequency-localization}. Since $\mu$ is real,
$\Re(T_\mu p_n)=T_\mu(\Re p_n)$, and therefore
\begin{align*}
        \|T_\mu p_n\|_1
        &\ge \|T_\mu(\Re p_n)\|_1 \\
        &\ge \|T_\mu h_n\|_1
             -\|\mu\|\,\|\Re p_n-h_n\|_1 .
\end{align*}
The error tends to zero along the above infinite sequence of scales. After discarding finitely many scales,
\begin{equation*}
        \|\mu\|\,\|\Re p_n-h_n\|_1
        \le \frac12c_0(\beta)c_\beta(\mu),
\end{equation*}
and hence
\begin{equation*}
        \|T_\mu p_n\|_1
        \ge \frac12c_0(\beta)c_\beta(\mu).
\end{equation*}
Since $\|p_n\|_1\le C_{\tau,\gamma}$ and
\begin{equation*}
        \operatorname{spec}p_n\subset[N_n,b_{\tau,\gamma}N_n^\tau],
        \qquad N_n\to\infty,
\end{equation*}
the result follows with $b=b_{\tau,\gamma}$ and
$C_{\rm loc}(\beta,\tau)=c_0(\beta)/(2C_{\tau,\gamma})$.
\end{proof}

\begin{cor}\label{cor:localized}
Let $\mu\in M_0(\T)$ be a real measure and suppose that
$\rho_\beta(\mu)>0$ for some $0<\beta<1$. Then, for every $\tau>1$, there exist $b=b(\tau)>0$ and an arbitrarily rapidly increasing sequence of positive integers $(k_j)$ such that, for
\begin{equation*}
        \Lambda=\bigcup_{j=1}^{\infty}[k_j,bk_j^\tau]\cap\mathbb Z,
\end{equation*}
the multiplier
$\widehat\mu|_\Lambda\in M(L^1_\Lambda(\T))$ has no natural spectrum.
\end{cor}
\begin{proof}
Choose numbers $0<c<\delta<\rho_\beta(\mu)$. By the definition of $\rho_\beta(\mu)$ there is an increasing sequence $m_j\to\infty$ such that
\begin{equation*}
        c_\beta(\mu^{*m_j})\ge\delta^{m_j}.
\end{equation*}
Since $\mu$ is Rajchman, there is $N(c)$ such that
$|\widehat\mu(n)|\le c$ for all $n\ge N(c)$. We choose $k_j$ inductively, arbitrarily rapidly increasing, so large that $k_j\ge N(c)$ and Proposition \ref{prop:localized}, applied to $\mu^{*m_j}$, gives
\begin{equation*}
\|T_{\mu^{*m_j}}:L^1_{[k_j,bk_j^\tau]}(\T)
\to L^1_{[k_j,bk_j^\tau]}(\T)\|
\ge C_{\rm loc}(\beta,\tau)c_\beta(\mu^{*m_j}).
\end{equation*}
Put $\Lambda=\bigcup_{j\ge1}[k_j,bk_j^\tau]\cap\mathbb Z$. Then
$\widehat\mu(\Lambda)\subset B(0,c)$ and
\begin{align*}
\|T_\mu^{m_j}:L^1_\Lambda(\T)\to L^1_\Lambda(\T)\|
&\ge C_{\rm loc}(\beta,\tau)c_\beta(\mu^{*m_j}) \\
&\ge C_{\rm loc}(\beta,\tau)\delta^{m_j}.
\end{align*}
Consequently,
\begin{equation*}
        r_{M(L^1_\Lambda(\T))}(T_\mu)\ge\delta>c,
\end{equation*}
whereas natural spectrum would have spectral radius at most $c$.
\end{proof}

\begin{cor}[Hereditary failure under arbitrary enlargements]\label{cor:upward}
For every $\tau>1$ there exist a positive Rajchman probability measure
$\mu$, a real Rajchman measure
\begin{equation*}
        \nu=\mu-m_\T,
\end{equation*}
of total mass zero, a constant $b=b(\tau)>0$, and an arbitrarily rapidly increasing sequence $(k_j)$ such that, with
\begin{equation*}
        \Lambda=\bigcup_{j=1}^{\infty}[k_j,bk_j^\tau]\cap\mathbb Z,
\end{equation*}
the multiplier $\widehat\nu|_\Gamma$ has no natural spectrum on
$L^1_\Gamma(\T)$ for every set $\Gamma\subset\mathbb Z$ satisfying
$\Lambda\subset\Gamma$.
\end{cor}
\begin{proof}
Let $\mu$ be the positive Rajchman probability measure carried by the set
$S_\infty$ described above. By Lemma \ref{lem:sumsets}, every convolution power $\mu^{*m}$ is carried by a set of lower Minkowski dimension zero. Hence
\begin{equation*}
        c_\beta(\mu^{*m})=1
\end{equation*}
for every $m\ge1$ and every $0<\beta<1$. Choose any sequence $m_j\to\infty$. Proposition \ref{prop:localized} allows us to choose an arbitrarily rapidly increasing sequence $(k_j)$ such that
\begin{equation*}
 \|T_\mu^{m_j}:L^1_{[k_j,bk_j^\tau]}(\T)
       \to L^1_{[k_j,bk_j^\tau]}(\T)\|
 \ge C_{\rm loc}(\beta,\tau).
\end{equation*}
For $\Lambda=\bigcup_j[k_j,bk_j^\tau]\cap\mathbb Z$ this yields
\begin{equation*}
        r_{M(L^1_\Lambda(\T))}(T_\mu)\ge1.
\end{equation*}
Since $\mu$ is a probability measure, $\|T_\mu\|\le1$, so the spectral radius is exactly one.

Let $m_\T$ be normalized Haar measure and put $\nu=\mu-m_\T$. Since
$0\notin\Lambda$, convolution with $m_\T$ vanishes on $L^1_\Lambda(\T)$; hence
\begin{equation*}
        T_\nu|_{L^1_\Lambda(\T)}=T_\mu|_{L^1_\Lambda(\T)}
\end{equation*}
and the spectral radius of $T_\nu$ on $L^1_\Lambda(\T)$ is one. If
$\Gamma\supset\Lambda$, then $L^1_\Lambda(\T)$ is an invariant subspace of
$L^1_\Gamma(\T)$, so
\begin{equation*}
        r_{M(L^1_\Gamma(\T))}(T_\nu)\ge1.
\end{equation*}
On the other hand,
\begin{equation*}
        \widehat\nu(0)=0,
        \qquad
        \widehat\nu(n)=\widehat\mu(n)\quad(n\ne0).
\end{equation*}
Moreover,
\begin{equation*}
        q:=\sup_{n\ne0}|\widehat\mu(n)|<1.
\end{equation*}
Indeed, equality at a non-zero frequency would force the corresponding character to be constant $\mu$-almost everywhere, contradicting the Rajchman property; the passage from pointwise strict inequality to the strict supremum uses
$\widehat\mu(n)\to0$. Thus
\begin{equation*}
        \overline{\widehat\nu(\Gamma)}\subset\{z\in\mathbb C:|z|\le q\},
\end{equation*}
whereas the spectrum of $T_\nu$ on $L^1_\Gamma(\T)$ has spectral radius at least one. Hence the spectrum is not natural.
\end{proof}

The analytic specialization of Corollary \ref{cor:multidim-zafran} and the preceding localized results motivate the following.
\begin{de}
We say that a set $\Gamma\subset\mathbb Z$ is a set of natural spectrum if every Rajchman multiplier from $M_0(L^1_\Gamma(\T))$ has a natural spectrum.
\end{de}
Corollary \ref{cor:upward} shows more than the failure of this property for one sparse union of intervals. For every $\tau>1$ one can prescribe an arbitrarily sparse sequence $(k_j)$ such that every set $\Gamma\subset\mathbb Z$ containing all the intervals
\begin{equation*}
        [k_j,bk_j^\tau]\cap\mathbb Z,
        \qquad j=1,2,\ldots,
\end{equation*}
fails the natural-spectrum property. The same zero-mass Rajchman measure is a witness for all such enlargements, including those which contain the zero frequency.

\section{Haar and Walsh--Stieltjes consequences of the mountain-river lemma}\label{sec:walsh}
The results in this section are one-dimensional and independent of the multiplier estimate. We apply the binary case of Lemma \ref{mrl} directly to Haar--Stieltjes blocks and then pass to Walsh--Stieltjes blocks by an elementary change of basis.

\subsection{Common notation and the Walsh--Haar transfer}
We identify the circle with the dyadic interval $[0,1)$, use the half-open dyadic intervals, and take $\N=\{1,2,\ldots\}$ in the Walsh notation. For classical background on Walsh--Stieltjes series, see \citep{Fine,SWS}. For a finite set $A\subset\N$ we write
\begin{equation*}
        w_A=\prod_{j\in A} r_j,
\end{equation*}
where $r_j$ is the $j$-th Rademacher function. If $\mu\in M([0,1))$, then
\begin{equation*}
        \widehat\mu(A)=\int_{[0,1)} w_A\,d\mu
\end{equation*}
denotes its Walsh--Stieltjes coefficient. We shall also use the full-carrier lower Minkowski dimension
\begin{equation*}
        \underline{\dim}_{M,\rm car}(\mu)=\inf\{\ldimM A: |\mu|([0,1)\setminus A)=0\}.
\end{equation*}
 This should not be confused with $c_\beta(\mu)$: the latter measures how much mass can be concentrated on a thin partial carrier, whereas $\underline{\dim}_{M,\rm car}(\mu)$ requires a set carrying the whole measure.

For a finite vector $a=(a_i)_{i=1}^{N}$ and an integer $m\ge1$ let
\begin{equation*}
        \|a\|_{W(m)}=\sum_{j=1}^{\min(m,N)} a_j^*,
\end{equation*}
where $(a_j^*)$ is the non-increasing rearrangement of $(|a_i|)$. Thus $\|a\|_{W(m)}$ is the sum of the $m$ largest absolute values. We shall use this notation for both Haar and Walsh blocks.

Let $\mathcal F_n$ be the family of dyadic intervals of length $2^{-n}$. If $I\in\mathcal F_n$ and $I_0,I_1$ are its left and right halves, put
\begin{equation*}
        \widetilde h_I=\chi_{I_0}-\chi_{I_1},
        \qquad
        c_I(\mu)=\int \widetilde h_I\,d\mu=\mu(I_0)-\mu(I_1).
\end{equation*}
The vector $(c_I(\mu))_{I\in\mathcal F_n}$ will be called the $n$-th Haar--Stieltjes block of $\mu$. If $\mu$ is a positive probability measure and
\begin{equation*}
        f_n=2^n\sum_{I\in\mathcal F_n}\mu(I)\chi_I
\end{equation*}
is the density of the conditional expectation of $\mu$ with respect to the dyadic sigma-algebra of level $n$, then the martingale difference $d_n=f_{n+1}-f_n$ satisfies
\begin{equation}\label{eq:haar-mart-diff}
        \|d_n\|_1=\sum_{I\in\mathcal F_n}|c_I(\mu)|.
\end{equation}
Indeed, on the two children of $I$ the function $d_n$ is equal to $\pm 2^n c_I(\mu)$.

For later use set
\begin{equation*}
        \mathcal W_n=\{A\subset\N:\max A=n+1\}.
\end{equation*}
This set has cardinality $2^n$. The following elementary change of basis will be used below.

\begin{lem}[Walsh--Haar transfer]\label{lem:walsh-haar-transfer}
For every $n$ and every finite measure $\mu\in M([0,1))$,
\begin{equation*}
        \bigl\|(c_I(\mu))_{I\in\mathcal F_n}\bigr\|_{W(m)}
        \le
        \bigl\|(\widehat\mu(A))_{A\in\mathcal W_n}\bigr\|_{W(m)}
\end{equation*}
for every integer $m\ge1$.
\end{lem}
\begin{proof}
Each function $\widetilde h_I$, $I\in\mathcal F_n$, belongs to the linear span of the Walsh functions $w_A$ with $\max A=n+1$. Thus, with the pairing $\langle f,g\rangle=\int_0^1 f(t)g(t)\,dt$,
\begin{equation*}
        c_I(\mu)=\sum_{A\in\mathcal W_n}\langle \widetilde h_I,w_A\rangle\widehat\mu(A).
\end{equation*}
The matrix $\bigl(\langle \widetilde h_I,w_A\rangle\bigr)_{I\in\mathcal F_n,A\in\mathcal W_n}$ has entries of absolute value $2^{-n}$, and the sum of absolute values in each row and each column is equal to $1$. Hence it is a contraction on both $\ell^1$ and $\ell^\infty$. The functional $\|\cdot\|_{W(m)}$ is an interpolation norm between $\ell^1$ and $\ell^\infty$; equivalently,
\begin{equation*}
        \|a\|_{W(m)}=\sup\left\{\left|\sum_i a_i b_i\right|:\|b\|_\infty\le1,\ \|b\|_1\le m\right\}.
\end{equation*}
Indeed, in the dual formula above the transpose matrix preserves both constraints $\|b\|_\infty\le1$ and $\|b\|_1\le m$, because it is also a contraction on $\ell^\infty$ and $\ell^1$. The desired inequality follows.
\end{proof}

\subsection{Haar and Walsh blocks}
This subsection contains the entire use of the combinatorial Lemma \ref{mrl} in the dyadic part of the paper. The lemma is applied directly only in the proof of Theorem \ref{thm:haar-block}. The Walsh statements which follow use it solely through that Haar theorem and Lemma \ref{lem:walsh-haar-transfer}.

\begin{thm}[Haar block lower estimate]\label{thm:haar-block}
Let $0<\beta<1/2$, and let $\mu$ be a non-zero finite positive measure on $[0,1)$ carried by a Borel set $A$ with $\ldimM A<\beta$. If
\begin{equation*}
        \frac{\beta}{1-\beta}<\lambda<1,
\end{equation*}
then there is a constant $C(\beta,\lambda)>0$ such that
\begin{equation*}
        \limsup_{n\to\infty}
        \bigl\|(c_I(\mu))_{I\in\mathcal F_n}\bigr\|_{W(\lfloor 2^{\lambda n}\rfloor)}
        \ge C(\beta,\lambda)\|\mu\|.
\end{equation*}
\end{thm}
\begin{proof}
Choose $\beta_0$ with $\ldimM A<\beta_0<\beta$. For infinitely many integers $k$, the set $A$ is covered by fewer than $2^{\beta_0 k}$ intervals of length $2^{-k}$. Replacing those intervals by the dyadic intervals of level $k$ which meet them increases the number only by an absolute factor. Since $\beta_0<\beta$, for all sufficiently large such $k$ we obtain a family $D_k\subset\mathcal F_k$ such that
\begin{equation*}
        A\subset\bigcup_{I\in D_k}I,
        \qquad
        \#D_k<2^{\beta k}.
\end{equation*}
Because $\mu$ is carried by $A$, the dyadic martingale of $\mu$ itself belongs to the class $M(\beta,k)$ at those levels.

We now use the combinatorial input of the paper. Since $\lambda>\beta/(1-\beta)$, choose $s\in(\beta,1)$ sufficiently close to $1$ that
\begin{equation*}
 \frac\beta\lambda<1-\frac\beta s,
\end{equation*}
and then fix
\begin{equation*}
 \frac\beta\lambda<\theta<1-\frac\beta s.
\end{equation*}
Choose $|\alpha|\in(s,1)$, put $\alpha=-|\alpha|$, and set $\rho=s/|\alpha|\in(0,1)$. Then $\alpha\rho+\beta=-s+\beta<0$. In the proof of Lemma \ref{mrl}, choose $\rho'<\rho$ sufficiently close to $\rho$ that
\begin{equation*}
 \theta<1-\frac{\beta}{|\alpha|\rho'}.
\end{equation*}
The binary case $d=1$ of that proof then gives a level $n>\theta k$ for which the turbulent vertices carry less than $\rho\|\mu\|$. No multiplier estimate is used in this argument.
 As in the proof of \eqref{wzn}, after choosing one of the two orientations we obtain a family $G_n\subset\mathcal F_n$ such that
\begin{equation*}
        \sum_{I\in G_n} |c_I(\mu)|
        \ge C_1(\beta,\lambda)\|\mu\|.
\end{equation*}
Every interval in $G_n$ has positive $\mu$-mass, hence it contains at least one interval from $D_k$ at level $k$. Therefore
\begin{equation*}
        \#G_n\le \#D_k<2^{\beta k}.
\end{equation*}
Because $\theta>\beta/\lambda$ and $n>\theta k$, for all sufficiently large $k$ we have $2^{\beta k}\le 2^{\lambda n}$. Hence the sum of the $\lfloor2^{\lambda n}\rfloor$ largest Haar coefficients at level $n$ is bounded from below by $C_1(\beta,\lambda)\|\mu\|$.
\end{proof}

After Theorem \ref{thm:haar-block}, the combinatorial lemma enters only through that theorem; the remaining passage is purely the Walsh--Haar change of basis.

\begin{thm}[Walsh block uncertainty principle]\label{thm:walsh-block}
Let $0<\beta<1/2$, and let $\mu$ be a non-zero finite positive measure on $[0,1)$ carried by a Borel set $A$ with $\ldimM A<\beta$. If
\begin{equation*}
        \frac{\beta}{1-\beta}<\lambda<1,
\end{equation*}
then there is a constant $C(\beta,\lambda)>0$ such that
\begin{equation*}
        \limsup_{n\to\infty}
        \bigl\|(\widehat\mu(A))_{A\in\mathcal W_n}\bigr\|_{W(\lfloor2^{\lambda n}\rfloor)}
        \ge C(\beta,\lambda)\|\mu\|.
\end{equation*}
\end{thm}
\begin{proof}
This is an immediate consequence of Theorem \ref{thm:haar-block} and Lemma \ref{lem:walsh-haar-transfer}.
\end{proof}

The preceding theorem can be read as a dimension criterion. It says that strong compressibility of all Walsh blocks is incompatible with concentration on sets of small lower Minkowski dimension.
\begin{cor}\label{cor:walsh-dim}
Let $0<\lambda<1$, and let $\mu$ be a non-zero positive finite measure on $[0,1)$. If
\begin{equation*}
        \lim_{n\to\infty}
        \bigl\|(\widehat\mu(A))_{A\in\mathcal W_n}\bigr\|_{W(\lfloor2^{\lambda n}\rfloor)}=0,
\end{equation*}
then
\begin{equation*}
        \underline{\dim}_{M,\rm car}(\mu)\ge \frac{\lambda}{1+\lambda}.
\end{equation*}
\end{cor}
\begin{proof}
If $\underline{\dim}_{M,\rm car}(\mu)<\lambda/(1+\lambda)$, then $\mu$ is carried by a Borel set $A$ with $\ldimM A<\lambda/(1+\lambda)$. Choose $\beta$ with
\begin{equation*}
        \ldimM A<\beta<\frac{\lambda}{1+\lambda}.
\end{equation*}
Then $\lambda>\beta/(1-\beta)$, contradicting Theorem \ref{thm:walsh-block}.
\end{proof}

\section{Two additional dyadic observations}\label{sec:further-dyadic}
This section records two independent observations leading to dimension estimates of a similar form. They are included not because they enter the proof of the multiplier theorem or of the mountain-river consequences, but because they lead to related Fourier-analytic statements on the dyadic group. The first starts from a classical entropy--predictability argument and turns it into estimates for Walsh--Stieltjes blocks. The second is an elementary support estimate for low-degree Walsh polynomials. Their common role in the paper is to complement the mountain-river method by showing, through two different mechanisms, that restrictions on dyadic Fourier data force lower bounds on the geometric size of every carrier.

\subsection{Entropy, prediction, and Haar martingale differences}
The entropy estimate below is a dyadic reformulation of classical entropy--error inequalities and is not claimed to be new in itself. For binary prediction, Hellman and Raviv \citep{HellmanRaviv} proved that the Bayes error is at most one half of the conditional entropy. Feder and Merhav \citep{FederMerhav} developed the corresponding sequential-prediction formulation and related average prediction error to the cardinality of a family of binary strings. Combined with the entropy-dimension viewpoint introduced by R\'enyi \citep{Renyi}, these facts yield the following Haar-martingale statement. We include its short proof because, after passing from Haar martingale differences to Walsh blocks, it produces Fourier-analytic consequences for Walsh--Stieltjes coefficients. In particular, it yields lower bounds for Ces\`aro averages of the $\ell^1$-norms of Walsh--Stieltjes blocks and implies that sufficiently sparse Walsh spectrum forces every carrier of a positive measure to have full lower Minkowski dimension. To the best of our knowledge, these Walsh--Stieltjes consequences have not previously been stated.

\begin{prop}[Entropy--predictability estimate for Haar martingale differences]\label{prop:haar-entropy}
Let $\mu$ be a positive probability measure on $[0,1)$. Then
\begin{equation}\label{eq:haar-entropy-estimate}
        \limsup_{N\to\infty}\frac1N\sum_{n=0}^{N-1}
        \sum_{I\in\mathcal F_n}|\mu(I_0)-\mu(I_1)|
        \ge
        1-\underline{\dim}_{M,\rm car}(\mu).
\end{equation}
Equivalently,
\begin{equation*}
        \limsup_{N\to\infty}\frac1N\sum_{n=0}^{N-1}\|f_{n+1}-f_n\|_1
        \ge
        1-\underline{\dim}_{M,\rm car}(\mu).
\end{equation*}
\end{prop}
\begin{proof}
For $I\in\mathcal F_n$ with $\mu(I)>0$ put
\begin{equation*}
        q_I=\frac{\mu(I_0)}{\mu(I)}.
\end{equation*}
Let
\begin{equation*}
        H_n(\mu)=-\sum_{I\in\mathcal F_n}\mu(I)\log_2\mu(I)
\end{equation*}
be the entropy of the dyadic projection of $\mu$, with the convention $0\log 0=0$. Then
\begin{equation*}
        H_{n+1}(\mu)-H_n(\mu)
        =\sum_{I\in\mathcal F_n}\mu(I)h(q_I),
\end{equation*}
where $h(q)=-q\log_2 q-(1-q)\log_2(1-q)$ is the binary entropy. On the other hand,
\begin{equation*}
        \sum_{I\in\mathcal F_n}|\mu(I_0)-\mu(I_1)|
        =\sum_{I\in\mathcal F_n}\mu(I)|2q_I-1|.
\end{equation*}
The classical binary entropy--error bound of Hellman and Raviv is, in the present notation,
\begin{equation*}
        \min(q,1-q)\le \frac12 h(q),
\end{equation*}
or equivalently
\begin{equation*}
        1-h(q)\le |2q-1|,\qquad 0\le q\le1.
\end{equation*}
For completeness, by symmetry it is enough to consider $q\in[1/2,1]$; the function $|2q-1|-(1-h(q))=2q-2+h(q)$ vanishes at $q=1/2$ and $q=1$, and is concave on this interval. Hence it is non-negative. Therefore
\begin{equation*}
        \sum_{I\in\mathcal F_n}|\mu(I_0)-\mu(I_1)|
        \ge 1-\bigl(H_{n+1}(\mu)-H_n(\mu)\bigr).
\end{equation*}
Summing from $n=0$ to $N-1$ yields
\begin{equation}\label{eq:entropy-telescope}
        \sum_{n=0}^{N-1}\sum_{I\in\mathcal F_n}|\mu(I_0)-\mu(I_1)|
        \ge N-H_N(\mu).
\end{equation}
Let $E$ be a Borel carrier of $\mu$, and let $N_d(E,2^{-N})$ be the number of dyadic intervals of length $2^{-N}$ which meet $E$. Since all mass of the dyadic projection of level $N$ is supported on those intervals, entropy is maximized when the masses are equal, and therefore
\begin{equation*}
        H_N(\mu)\le \log_2 N_d(E,2^{-N}).
\end{equation*}
The dyadic covering numbers give the same lower Minkowski dimension. Hence, along a sequence of $N$'s realizing the lower Minkowski dimension of $E$, \eqref{eq:entropy-telescope} gives
\begin{equation*}
        \limsup_{N\to\infty}\frac1N\sum_{n=0}^{N-1}\sum_{I\in\mathcal F_n}|\mu(I_0)-\mu(I_1)|
        \ge 1-\ldimM E.
\end{equation*}
Taking carriers $E$ with $\ldimM E$ arbitrarily close to $\underline{\dim}_{M,\rm car}(\mu)$ proves the proposition.
\end{proof}

\begin{rem}[On positivity]\label{rem:haar-positivity}
The positivity assumption in Proposition \ref{prop:haar-entropy} is essential for the entropy argument. For signed measures, dyadic cancellations may make the raw Haar martingale differences small even when the total variation is carried by a very thin set. The real-measure estimate in Theorem \ref{dwa} avoids this difficulty by first isolating a one-signed part of the measure and then using adapted Hardy atoms rather than the Haar martingale differences of the signed measure itself.
\end{rem}

\begin{rem}\label{rem:distribution-function-haar}
If $\mu=dF$, where $F(t)=\mu([0,t))$, then for $I=[a_I,b_I)$ and $m_I=(a_I+b_I)/2$ we have
\begin{equation*}
        \mu(I_0)-\mu(I_1)=2F(m_I)-F(a_I)-F(b_I).
\end{equation*}
Thus Proposition \ref{prop:haar-entropy} says that a distribution function carried by a thin set has large average dyadic second differences. In particular, if $\underline{\dim}_{M,\rm car}(\mu)=0$, then the limsup in \eqref{eq:haar-entropy-estimate} is equal to $1$, the maximal possible value.
\end{rem}

We now extract the Fourier-analytic consequence for which the preceding entropy formulation was introduced. Although the entropy--predictability estimate itself is classical, its transfer to Walsh--Stieltjes coefficients gives the following dimension bound for averaged Walsh blocks.
\begin{cor}[Averaged Walsh variation and carrier dimension]\label{cor:walsh-entropy}
Let $\mu$ be a non-zero finite positive measure on $[0,1)$. Then
\begin{equation*}
        \limsup_{N\to\infty}\frac1N\sum_{n=0}^{N-1}
        \bigl\|(\widehat\mu(A))_{A\in\mathcal W_n}\bigr\|_{\ell^1}
        \ge
        \bigl(1-\underline{\dim}_{M,\rm car}(\mu)\bigr)\|\mu\|.
\end{equation*}
\end{cor}
\begin{proof}
Normalize $\mu$ to have total mass one, use Lemma \ref{lem:walsh-haar-transfer} with $m=2^n$, and apply Proposition \ref{prop:haar-entropy}. The stated form follows by homogeneity.
\end{proof}

\paragraph{Comparison with Theorem \ref{thm:walsh-block}.}
Corollary \ref{cor:walsh-entropy} and Theorem \ref{thm:walsh-block} both turn geometric thinness of a carrier into a lower bound for Walsh--Stieltjes blocks, but they measure different features of those blocks and neither result implies the other. Theorem \ref{thm:walsh-block} is a sparse, scale-by-scale statement: if the carrier has lower Minkowski dimension below $\beta<1/2$, then on infinitely many levels the sum of only the $2^{\lambda n}$ largest coefficients, with $\lambda>\beta/(1-\beta)$, already captures a fixed proportion of $\|\mu\|$. It therefore detects concentration, or compressibility, inside individual Walsh blocks. It does not, however, control how frequently the good levels occur.

By contrast, Corollary \ref{cor:walsh-entropy} applies throughout the full range $0\le \underline{\dim}_{M,\rm car}(\mu)\le1$ and gives information averaged over all levels. Its conclusion concerns the full $\ell^1$-norm of each block and therefore permits that norm to be spread among all $2^n$ coefficients; it gives no sparse-subblock conclusion. Thus the mountain-river argument is stronger with respect to concentration inside selected blocks, whereas the entropy argument is stronger with respect to the range of dimensions and the density of scales. Their complementary nature is the reason for including the entropy-based result separately.

For sparse Walsh spectra the averaged form is often the most convenient one.
Let $W$ be a family of finite subsets of $\N$. We write $M_W([0,1))$ for the space of all measures $\mu\in M([0,1))$ such that
\begin{equation*}
        \widehat\mu(A)=0\qquad\text{for all }A\notin W.
\end{equation*}
\begin{cor}[Sparse Walsh spectra]\label{cor:sparse-walsh}
Assume that
\begin{equation*}
        \#\{A\in W:A\subset\{1,\ldots,k\}\}=o(k).
\end{equation*}
Then every non-zero positive measure $\mu\in M_W([0,1))$ satisfies
\begin{equation*}
        \underline{\dim}_{M,\rm car}(\mu)=1.
\end{equation*}
\end{cor}
\begin{proof}
Normalize $\mu$ to have total mass $1$. Since $\mu\in M_W([0,1))$ and $|\widehat\mu(A)|\le1$, we have
\begin{equation*}
        \frac1N\sum_{n=0}^{N-1}
        \bigl\|(\widehat\mu(A))_{A\in\mathcal W_n}\bigr\|_{\ell^1}
        \le
        \frac{\#\{A\in W:A\subset\{1,\ldots,N\}\}}{N}
        \longrightarrow0.
\end{equation*}
Corollary \ref{cor:walsh-entropy} therefore gives
\begin{equation*}
        1-\underline{\dim}_{M,\rm car}(\mu)\le0.
\end{equation*}
Thus $\underline{\dim}_{M,\rm car}(\mu)=1$.
\end{proof}

The last corollary gives a simple family of singular measures with sparse Walsh spectrum but full lower Minkowski dimension.
\begin{cor}[Sparse Walsh--Riesz products]\label{cor:walsh-riesz-products}
Let $(n_j)$ be an increasing sequence of positive integers such that
\begin{equation*}
        2^j=o(n_j).
\end{equation*}
Let
\begin{equation*}
        d\mu=\prod_{j=1}^{\infty}(1+a_j r_{n_j})\,dt,
        \qquad -1\le a_j\le1,
\end{equation*}
be the corresponding Walsh--Riesz product. Then
\begin{equation*}
        \underline{\dim}_{M,\rm car}(\mu)=1.
\end{equation*}
In particular, if $a_j\to0$ and $\sum_j |a_j|^2=\infty$, then $\mu$ is a singular Walsh--Rajchman measure whose every Borel carrier has full lower Minkowski dimension.
\end{cor}
\begin{proof}
The Walsh spectrum of the product is contained in
\begin{equation*}
        W=\{A\subset\{n_1,n_2,\ldots\}:\#A<\infty\}.
\end{equation*}
If $n_j\le N<n_{j+1}$, then
\begin{equation*}
 \#\{A\in W:A\subset\{1,\ldots,N\}\}=2^j,
\end{equation*}
and therefore
\begin{equation*}
 \frac{\#\{A\in W:A\subset\{1,\ldots,N\}\}}{N}
 \le\frac{2^j}{n_j}\longrightarrow0.
\end{equation*}
Thus the hypothesis of Corollary \ref{cor:sparse-walsh} holds for all $N$, not merely along the distinguished subsequence, and that corollary gives $\underline{\dim}_{M,\rm car}(\mu)=1$. The final assertion is the standard Kakutani dichotomy for Walsh--Riesz products, see for instance \citep[Chapter 7]{grmc}: $a_j\to0$ is exactly the Rajchman condition for the Walsh coefficients, while $\sum_j |a_j|^2=\infty$ implies singularity with respect to Lebesgue measure.
\end{proof}

\subsection{Restricted Walsh degree}
This argument uses only a support bound for low-degree Walsh polynomials.

A different and very transparent form of Walsh sparsity is obtained by restricting the degree of the Walsh products which are allowed to occur up to a given level. We shall use the following elementary fact.

\begin{lem}[Support of a low-degree Walsh polynomial]\label{lem:degree-support}
Let $P$ be a non-zero Walsh polynomial on the dyadic group $
\{0,1\}^{n}$ of degree at most $d$, $0\le d\le n$. Then $P$ is non-zero on at least $2^{n-d}$ points.
\end{lem}
\begin{proof}
We argue by induction on $n$. The assertion is clear for $n=0$, and it is trivial when $d=n$. Assume therefore that $d<n$ and write
\[
        P(x_1,\ldots,x_n)=P_0(x_1,\ldots,x_{n-1})+r_n P_1(x_1,\ldots,x_{n-1}),
\]
where $\deg P_0\le d$ and $\deg P_1\le d-1$. If both restrictions $P_0+P_1$ and $P_0-P_1$ are non-zero, then the induction hypothesis applies with degree at most $d\le n-1$ to each restriction. Each is non-zero on at least $2^{n-1-d}$ points, so $P$ is non-zero on at least $2^{n-d}$ points. If one restriction is identically zero, the other is twice a non-zero polynomial of degree at most $d-1$ and is supported over one half of the cube. The induction hypothesis gives at least $2^{(n-1)-(d-1)}=2^{n-d}$ non-zero points.
\end{proof}

\begin{prop}[Degree-restricted Walsh spectra]\label{prop:degree-restricted}
Let $F:\N\to\N$ be non-decreasing, and let $\mu$ be a non-zero finite real or complex measure on $[0,1)$. Suppose that
\begin{equation*}
        \widehat\mu(A)=0
        \qquad
        \text{whenever } A\subset\{1,\ldots,n\}
        \text{ and } |A|>F(n).
\end{equation*}
Then
\begin{equation*}
        \underline{\dim}_{M,\rm car}(\mu)
        \ge
        1-\limsup_{n\to\infty}\frac{F(n)}{n}.
\end{equation*}
In particular, if $F(n)=o(n)$, then every Borel carrier of $\mu$ has full lower Minkowski dimension.
\end{prop}
\begin{proof}
Let $f_n$ be the density of the conditional expectation of $\mu$ with respect to the dyadic sigma-algebra generated by $r_1,\ldots,r_n$. Equivalently,
\begin{equation*}
        f_n=\sum_{A\subset\{1,\ldots,n\}}\widehat\mu(A)w_A.
\end{equation*}
By the assumption, $f_n$ is a Walsh polynomial of degree at most $F(n)$. Since $\mu$ is non-zero, $f_n$ is not identically zero for all sufficiently large $n$. Indeed, otherwise $\mu$ would vanish on all dyadic intervals from an infinite sequence of levels, hence on the algebra generated by dyadic intervals. Lemma \ref{lem:degree-support}, applied with $d=\min(F(n),n)$ and with such large $n$, implies that $f_n$ is non-zero on at least
\begin{equation*}
        2^{n-\min(F(n),n)}
\end{equation*}
dyadic atoms of level $n$.

Let $E$ be a Borel carrier of $\mu$. If $f_n$ is non-zero on a dyadic atom $I$, then $\mu(I)\ne0$. Hence $|\mu|(I)>0$, and every Borel carrier $E$ of $\mu$ must satisfy $E\cap I\ne\emptyset$. Let $N_{\rm dyad}(E,2^{-n})$ be the number of level-$n$ dyadic intervals which meet $E$. Then
\begin{equation*}
        N_{\rm dyad}(E,2^{-n})\ge 2^{n-\min(F(n),n)}.
\end{equation*}
An arbitrary interval of length $2^{-n}$ meets at most three level-$n$ dyadic intervals (the harmless third interval can occur only through endpoint conventions). Therefore the covering number from the definition of lower Minkowski dimension satisfies
\begin{equation*}
 N(E,2^{-n})\ge\frac13N_{\rm dyad}(E,2^{-n}),
\end{equation*}
and consequently
\begin{equation*}
        \ldimM E
        \ge \liminf_{n\to\infty}\frac{n-\min(F(n),n)}{n}
        = 1-\limsup_{n\to\infty}\min\left(\frac{F(n)}{n},1\right).
\end{equation*}
This implies the stated estimate, which is non-trivial exactly when $\limsup F(n)/n\le1$.
\end{proof}

\begin{rem}[Sharpness]\label{rem:degree-sharpness}
The estimate in Proposition \ref{prop:degree-restricted} is sharp. Let $S\subset\N$ have asymptotic density $\alpha$, and let $\mu$ be the Haar probability measure on the dyadic subcube on which the coordinates in $S$ are fixed, say $r_j=1$ for $j\in S$. Then
\begin{equation*}
        \widehat\mu(A)=0\qquad\text{unless }A\subset S.
\end{equation*}
Thus the hypothesis of Proposition \ref{prop:degree-restricted} holds with
\begin{equation*}
        F(n)=\#(S\cap\{1,\ldots,n\})\sim \alpha n.
\end{equation*}
At level $n$ the carrier meets exactly $2^{n-F(n)}$ dyadic atoms, and therefore its lower Minkowski dimension is $1-\alpha$. Conversely, Proposition \ref{prop:degree-restricted} applied to this choice of $F$ shows that every full carrier has lower Minkowski dimension at least $1-\alpha$; hence the displayed dyadic subcube attains the bound exactly.
\end{rem}

\section{Open problems}
The first two questions are related to the classical Riesz-product example proved in Appendix \ref{app:riesz}. That example is not used elsewhere in the paper.
\begin{enumerate}
    \item What is the spectrum of the classical Riesz product from Theorem \ref{trz}, treated as a multiplier on $H^{1}(\T)$? Its spectrum in $M(\T)$ is the whole unit disc; see \citep[Chapter 6]{grmc} on independent power measures.
    \item Can one find a multiplier on $H^{1}(\T)$ with non-natural spectrum within the class of Riesz products in $M_{0}(\T)$?
    \item Can one find a Rajchman probability measure whose spectrum as an operator on $H^1_0(\T)$ is the unit circle or the unit disc?
    \item Can the lower Minkowski dimension in Theorem \ref{dwa} be replaced by Hausdorff dimension?
    \item How far from optimal is the norm estimate in Theorem \ref{dwa}? In particular, how does the best constant behave as $\beta\uparrow d$?
    \item Can the localized intervals in Corollary \ref{cor:upward} be replaced by intervals $[N,bN]$ with a fixed $b>1$, or more generally by substantially shorter intervals? The present argument yields $[N,bN^\tau]$ for every $\tau>1$ but does not reach the linear endpoint.
\end{enumerate}

\appendix
\section{A classical Riesz-product example}\label{app:riesz}
In this appendix only, $H^1(\T)$ denotes the analytic Hardy space including the constants. Once the Rajchman assumption is dropped, non-natural spectrum for an $H^1(\T)$ multiplier can already be obtained from the standard Riesz product; compare Graham's classical Riesz-product proof of the Wiener--Pitt theorem \citep{graham}. The functional-calculus step is classical; we include the short counting argument needed to connect it rigorously with Klemes's characterization of idempotent $H^1$ multipliers.

\begin{thm}\label{trz}
Let
\begin{equation*}
    \mu=\mathop{\ast\!\!\operatorname{-weak}}\lim_{n\to\infty}
    \prod_{k=1}^{n}\left(1+\cos(3^{k}t)\right)
\end{equation*}
be the standard Riesz product. Then
\begin{equation*}
 \sigma_{B(H^{1}(\T))}(T_{\mu})
 \neq \overline{\widehat{\mu}(\mathbb Z_{\ge0})}.
\end{equation*}
\end{thm}

\begin{proof}
Assume to the contrary that
\begin{equation*}
 \sigma_{B(H^{1}(\T))}(T_{\mu})
 =\overline{\widehat{\mu}(\mathbb Z_{\ge0})}
 =\{0\}\cup\left\{2^{-r}:r\ge1\right\}\cup\{1\}.
\end{equation*}
Separate the point $1/4$ from the rest of this compact set by disjoint open sets and apply the holomorphic functional calculus in $B(H^1(\T))$. The resulting idempotent still commutes with rotations and hence is a Fourier multiplier. Its symbol is the indicator of
\begin{equation*}
 E=\{n\in\mathbb Z_{\ge0}:\widehat\mu(n)=1/4\}.
\end{equation*}
The balanced ternary expansion is unique, and the Fourier expansion of the Riesz product therefore gives
\begin{equation}\label{eq:riesz-level-two}
 E=\bigl\{n\in\mathbb Z_{\ge0}:n=\pm3^j\pm3^k
 \text{ for distinct }j,k\in\mathbb N\bigr\}.
\end{equation}

By Theorem 1.2 of Klemes \citep{kle}, the support of an idempotent multiplier on $H^1(\T)$ is a finite Boolean combination of finite sets, arithmetic progressions, and lacunary sequences.  Any finite Boolean combination of this form differs from an eventually periodic set $P_0\subset\mathbb Z_{\ge0}$ by a subset of a finite set together with a finite union of lacunary sequences. Indeed, replace every finite or lacunary generator in a Boolean formula by the empty set; the value of the formula can change only at a point belonging to one of those generators. Such an exceptional set contains $O(\log N)$ integers in $[1,N]$.

The set $E$ has density zero. Indeed, if $0<\pm3^j\pm3^k\le N$ with $j\ne k$, then $\max(j,k)\le C\log N$; hence \eqref{eq:riesz-level-two} gives $\#(E\cap[1,N])=O((\log N)^2)$.  Consequently the periodic tail of $P_0$ must be empty, since every non-empty periodic tail has positive density whereas the symmetric difference $E\triangle P_0$ has density zero. Thus $P_0$ is finite and can be absorbed into the exceptional finite set. Klemes's characterization would therefore imply
\begin{equation}\label{eq:riesz-upper-count}
 \#(E\cap[1,N])=O(\log N).
\end{equation}
On the other hand, for every $m\ge2$ the integers
\begin{equation*}
 3^j+3^k,\qquad 1\le k<j\le m,
\end{equation*}
are pairwise distinct, belong to $E$, and do not exceed $2\cdot3^m$. Thus
\begin{equation*}
 \#\bigl(E\cap[1,2\cdot3^m]\bigr)\ge\binom m2,
\end{equation*}
which contradicts \eqref{eq:riesz-upper-count}. Therefore the spectrum of $T_\mu$ is not natural.
\end{proof}

\section*{Acknowledgments}
The authors thank Fedor Nazarov for valuable discussion during the work on this paper.

\section*{Statements and Declarations}

\noindent\textbf{Funding.}
\begin{sloppypar}
The research of Michał Wojciechowski was partially supported by the National Science Centre, Poland, and Austrian Science Foundation FWF joint CEUS programme. National Science Centre project no.~2020/02/Y/ST1/00072 and FWF project no.~I5231.
\end{sloppypar}

\noindent\textbf{Competing interests.}
The authors have no relevant financial or non-financial interests to disclose.

\noindent\textbf{Author contributions.}
All authors contributed to the mathematical development and preparation of the manuscript. All authors read and approved the final manuscript.

\noindent\textbf{Data availability.}
No datasets were generated or analysed during the current study.

\noindent\textbf{Use of artificial intelligence.}
During the preparation of this work, the authors used ChatGPT (OpenAI) to assist with language editing and manuscript preparation. All AI-assisted material was reviewed and, where necessary, revised by the authors, who take full responsibility for the content of the paper.


\begin{thebibliography}{99}

\bibitem{BS}
Barral, J., Seuret, S. (eds.): Recent Developments in Fractals and Related Fields. Applied and Numerical Harmonic Analysis. Birkh\"auser, Boston (2010)

\bibitem{B}
Bluhm, C.: Liouville numbers, Rajchman measures, and small Cantor sets. Proc. Am. Math. Soc. 128(9), 2637--2640 (2000). \url{https://doi.org/10.1090/S0002-9939-00-05452-6}

\bibitem{FederMerhav}
Feder, M., Merhav, N.: Relations between entropy and error probability. IEEE Trans. Inf. Theory 40(1), 259--266 (1994). \url{https://doi.org/10.1109/18.272494}

\bibitem{FeffermanStein}
Fefferman, C., Stein, E.M.: $H^p$ spaces of several variables. Acta Math. 129, 137--193 (1972). \url{https://doi.org/10.1007/BF02392215}

\bibitem{Fine}
Fine, N.J.: Fourier--Stieltjes series of Walsh functions. Trans. Am. Math. Soc. 86, 246--255 (1957). \url{https://doi.org/10.1090/S0002-9947-1957-0087645-5}

\bibitem{graham}
Graham, C.C.: A Riesz product proof of the Wiener--Pitt theorem. Proc. Am. Math. Soc. 44(2), 312--314 (1974). \url{https://doi.org/10.1090/S0002-9939-1974-0340972-1}

\bibitem{grmc}
Graham, C.C., McGehee, O.C.: Essays in Commutative Harmonic Analysis. Springer, New York (1979)

\bibitem{HellmanRaviv}
Hellman, M.E., Raviv, J.: Probability of error, equivocation, and the Chernoff bound. IEEE Trans. Inf. Theory 16(4), 368--372 (1970). \url{https://doi.org/10.1109/TIT.1970.1054466}

\bibitem{kle}
Klemes, I.: Idempotent multipliers of $H^1(\T)$. Can. J. Math. 39(5), 1223--1234 (1987). \url{https://doi.org/10.4153/CJM-1987-062-5}

\bibitem{Renyi}
R\'enyi, A.: On the dimension and entropy of probability distributions. Acta Math. Acad. Sci. Hungar. 10, 193--215 (1959). \url{https://doi.org/10.1007/BF02063299}

\bibitem{SWS}
Schipp, F., Wade, W.R., Simon, P., P\'al, J.: Walsh Series: An Introduction to Dyadic Harmonic Analysis. Adam Hilger, Bristol (1990)

\bibitem{SleddStegenga}
Sledd, W.T., Stegenga, D.A.: An $H^1$ multiplier theorem. Ark. Mat. 19, 265--270 (1981). \url{https://doi.org/10.1007/BF02384484}

\bibitem{Stein}
Stein, E.M.: Harmonic Analysis: Real-Variable Methods, Orthogonality, and Oscillatory Integrals. Princeton Mathematical Series, vol. 43. Princeton University Press, Princeton (1993)

\bibitem{wp}
Wiener, N., Pitt, H.R.: On absolutely convergent Fourier--Stieltjes transforms. Duke Math. J. 4(2), 420--436 (1938). \url{https://doi.org/10.1215/S0012-7094-38-00435-1}

\bibitem{Zafran}
Zafran, M.: On the spectra of multipliers. Pac. J. Math. 47(2), 609--626 (1973). \url{https://doi.org/10.2140/pjm.1973.47.609}

\bibitem{Zaf2}
Zafran, M.: Measures as convolution operators on $H^1$ and $\mathrm{Lip}\,\alpha$. J. Funct. Anal. 29(2), 160--183 (1978). \url{https://doi.org/10.1016/0022-1236(78)90004-6}

\bibitem{Zygmund}
Zygmund, A.: Trigonometric Series, vols. I--II, 3rd edn. Cambridge University Press, Cambridge (2002)

\end{thebibliography}
\end{document}